\documentclass[12pt]{amsart}

\usepackage{amsthm}

\theoremstyle{plain}
\newtheorem{theorem}{Theorem}[section]

\newtheorem{proposition}[theorem]{Proposition}
\newtheorem{corollary}[theorem]{Corollary}

\theoremstyle{definition}
\newtheorem{example}[theorem]{Example}
\newtheorem{remark}[theorem]{Remark}
\newtheorem{definition}[theorem]{Definition} 

\usepackage[matrix,arrow,curve,cmtip]{xy}
  \usepackage{amsfonts}
   \usepackage{amsmath}
  \usepackage{amssymb}
 \numberwithin{equation}{section}
  \usepackage{hyperref}
\usepackage{enumerate,xspace}

\renewcommand{\phi}{\varphi}
\renewcommand{\epsilon}{\varepsilon}

\newcommand{\cb}{\ensuremath{\mathcal B}\xspace}
\newcommand{\cc}{\ensuremath{\mathcal C}\xspace}
\newcommand{\cd}{\ensuremath{\mathcal D}\xspace}

\newcommand{\ck}{\ensuremath{\mathcal K}\xspace}

\newcommand{\ba}{\ensuremath{\mathbb A}\xspace}
\newcommand{\bd}{\ensuremath{\mathbb D}\xspace}

\newcommand{\Span}{\ensuremath{\mathbf{Span}}\xspace}

\newcommand{\Cat}{\ensuremath{\mathbf{Cat}}\xspace}
\newcommand{\Vect}{\ensuremath{\mathbf{Vect}}\xspace}

\DeclareMathOperator{\ob}{ob}

\def\ox{\otimes}
\def\x{\times}

\def\1c#1{\stackrel{#1}{\to}}

\begin{document}


 \title{On monads and warpings}

\author{Stephen Lack}
\address{Department of Mathematics\\ Macquarie University}
\email{steve.lack@mq.edu.au}

\author{Ross Street}
\address{Department of Mathematics\\ Macquarie University}
\email{ross.street@mq.edu.au}

\begin{abstract}
We explain the sense in which a warping on a monoidal category is the same as a pseudomonad on the corresponding one-object bicategory, and we describe extensions of this to the setting of skew monoidal categories: these are a generalization of monoidal categories in which the associativity and unit maps are not required to be invertible. Our analysis leads us to describe a normalization process for skew monoidal categories, which produces a universal skew monoidal category for which the right unit map is invertible. 
\end{abstract}


\maketitle

\section{Introduction}

If \cc is a monoidal category with tensor product $\ox$, and $T\colon\cc\to\cc$ is a functor, then one can define a new product $\boxtimes$ on \cc via the formula 
$$A\boxtimes B = TA\ox B.$$
In order for this to define a new monoidal structure on \cc, further structure on \cc is required. The notion of {\em warping}, introduced in \cite{BookerStreet-Tannaka}, is designed to do just that: if $T$ is a warping then \cc becomes monoidal via the ``warped'' tensor product $\boxtimes$ defined above. 

While the notion of warping is quite restrictive, the {\em skew warpings} of \cite{skew} are far more common: for example, if $T$ has a monad structure, and this monad is opmonoidal \cite{Moerdijk-Hopf-monads,McCrudden-opmonoidal}, in the sense that there are suitably coherent maps $T(A\ox B)\to TA\ox TB$ and $TI\to I$, then $T$ is a skew warping. In particular, if $H$ is a bialgebra, then the functor $H\ox-\colon\Vect\to\Vect$ has a skew warping structure. 

The price of this extra generality is that the warped tensor product no long gives a monoidal structure, but only a {\em skew monoidal} one, in the sense of \cite{skew} (called {\em left skew monoidal} in \cite{Szlachanyi-skew}). These skew monoidal categories are similar to monoidal categories, except that the associativity and unit structure morphisms are not required to be invertible. The key insight of \cite{Szlachanyi-skew} is that these skew monoidal categories can be used to provide a valuable new characterization of bialgebroids; this was extended in \cite{skew} to the case of quantum categories.

We have been studying skew monoidal categories in a series of papers \cite{skew,Street-SkewClosed,skewEH,skewcoherence}, but have so far only scratched the surface of this remarkable theory, which seems to stem from the fact that skew monoidal categories are at the same time a generalization of monoidal categories and of categories.  

While skew warpings and skew monoidal structures are quite recent, monads have of course been a central topic in category theory for decades, and have been generalized in many directions. For example, monads can be defined in any bicategory \cite{bicategories}, and while monads in \Cat are just ordinary monads, monads in \Span are categories. Generalizing in a different direction, one can consider monads not just on categories but on 2-categories or bicategories, and in this context one often has weaker structures called pseudomonads; still more generally, there are various lax notions of monad. 

Most of these generalizations rely, directly or indirectly, on the fact that (ordinary) monads are the same as monoids in a monoidal category of endo\-functors. But there is also another approach, which  has largely been developed and promoted by Manes, for example in \cite{Manes-AlgebraicTheories}; but see also Walters' thesis \cite{Walters-thesis}. In this approach, one does not specify a functor at all; rather, for each object $A$ of the category \cc one gives an object $DA$ and a morphism $K_A\colon A\to DA$, and for each morphism $f\colon A\to DB$ one gives a morphism $Tf\colon DA\to DB$. 

One feature of this approach is that, whereas the usual definition of monad involves an associative multiplication $D\circ D\to D$ and so requires the formation of  $D\circ D$ and $D\circ D\circ D$, in Manes' approach, these iterates of $D$ are not needed. Thus Marmolejo and Wood use the epithet ``no iteration'' to refer to this approach to monads, when in \cite{MarmolejoWood-NoIterationPsm} they modify the theory to deal with pseudomonads. Since this is a little unwieldy, we shall replace ``no iteration'' by ``mw-''. We leave to the reader the question of whether these letters denote Manes and Walters, Marmolejo and Wood, or something else entirely.

The goal of this paper is to describe a close relationship between warpings and skew warpings on the one hand, and mw-monads and pseudo-mw-monads on the other.

Perhaps the simplest result to state is this:

\begin{quote}
 Let \cc be a monoidal category, and $\Sigma\cc$ the corresponding one-object bicategory. A warping on \cc is the same as a pseudomonad on $\Sigma\cc$.  
\end{quote}

\noindent
We prove this in Corollary~\ref{cor:main} below.
We  could equally have put pseudo-mw-monad rather than pseudomonad since, as proved in  \cite{MarmolejoWood-NoIterationPsm}, these amount to the same thing. 

This correspondence between pseudo-mw-monads and pseudomonads depends heavily on the invertibility of certain structure maps. If one weakens this requirement, the resulting notion of skew mw-monad is no longer equivalent to any lax version of ordinary pseudomonads. Nonetheless these skew mw-monads seem to be an interesting structure:

\begin{quote}
 Let \cc be a monoidal category, and $\Sigma\cc$ the corresponding one-object bicategory. A skew warping on \cc is the same as a skew mw-monad on $\Sigma\cc$.  
\end{quote}

These connections between (possibly skew) warpings and (higher) mw-monads shed light on both. In one direction, it shows that the ``warped'' monoidal structure involving $\boxtimes$ is really a sort of Kleisli construction, it suggests that one should consider ``algebras'' for skew warpings, and it suggests that warpings should be considered on bicategories as well as monoidal categories. In the other, it makes clear that some of the axioms for mw-pseudomonads are redundant, and suggests considering lax/skew variants as well. 

In the final section of the paper, we describe a universal process whereby a skew monoidal category can be replaced by one which is right normal, in the sense that the right unit constraint is invertible. We call this process (right) normalization, and we use it to give a formal account of the relationship between monads and mw-monads. 

\section{Review of mw-monads} \label{sect:monads}

In this section, we recall the definition of mw-monad, and its relationship to ordinary monads. 

The usual notion of monad on a category \cc consists of a functor $D\colon\cc\to\cc$ equipped with natural transformations $m\colon D^2\to D$ and $K\colon 1\to D$ satisfying associativity and unit laws. 

\begin{definition}
An {\em mw-monad} on \cc, consists of the following structure:
\begin{itemize}
\item a function $D\colon\ob\cc\to\ob\cc$
\item functions $T\colon\cc(X,DY)\to\cc(DX,DY)$ assigning to each morphism $f\colon X\to DY$ a morphism $Tf\colon DX\to DY$
\item a morphism $K=K_X\colon X\to DX$ for each $X$
\item subject to the following equations:
  \begin{align*}
    Tg\circ Tf &= T(Tg\circ f) \\
    Tf\circ K &= f \\
    TK_X &= 1_{DX}.
  \end{align*}
\end{itemize}  
\end{definition}

This determines a monad on \cc as follows. The endofunctor is defined on objects using $D$, and sends a morphism $f\colon X\to Y$ to $T(K_Y\circ f)$. The components of the unit are given by the $K_X$. The component at $X$ of the multiplication is $T(1_{DX})$. Conversely, for any monad $D$ on \cc with multiplication $M$ and unit $K$, we get an mw-monad by defining $Tf\colon DX\to DY$ to be $Df\colon DX\to D^2Y$ composed with the multiplication $D^2Y\to DY$. These constructions are mutually inverse: see \cite{Manes-AlgebraicTheories}.

These mw-monads are in some sense more closely related to their Kleisli categories than in the usual approach. Given an mw-monad as above, the Kleisli category $\cc_T$ has the same objects as \cc, with $\cc_T(X,Y)=\cc(X,DY)$; the identity on $X$ is $K_X$, while the composite of $f\colon X\to DY$ and $g\colon Y\to DZ$ is $Tg\circ f$.

It is also possible to reformulate the usual notion of algebra for a monad in terms of the mw-monad. This is done in the following definition.

\begin{definition}
Given an mw-monad as above, an {\em algebra} consists of an object $A$, together with  functions $E\colon\cc(X,A)\to\cc(DX,A)$ such that, for all $g\colon Y\to A$ and $f\colon X\to DY$, we have $Eg\circ K_Y=g$ and $Eg\circ Tf=E(Eg\circ f)$.
\end{definition}

\section{Skew bicategories}

There is an evident common generalization of the notions of bicategory and skew monoidal category, which we shall tentatively call a {\em skew bicategory}, although there are also richer structures which may deserve this name. At this stage, the only motivation for the definition is to have a common setting in which to discuss bicategories and skew monoidal categories. In any case, for this paper, a skew bicategory consists of:
\begin{itemize}
\item objects $X,Y,Z,\ldots$
\item hom-categories $\cb(X,Y)$ for all objects $X$ and $Y$
\item functors $M\colon\cb(Y,Z)\x\cb(X,Y)\to\cb(X,Z)$ 
\item functors $j\colon 1\to \cb(X,X)$
\item (not necessarily invertible) natural transformations
$$\xymatrix{
\cb(Y,Z)\x \cb(X,Y)\x\cb(W,X) \ar[r]^-{M\x1}_{~}="1" \ar[d]_{1\x M} & 
\cb(X,Z)\x\cb(W,X) \ar[d]^{M} \\
\cb(Y,Z)\x\cb(W,Y) \ar[r]_{M}^{~}="2" & \cb(W,Z)
\ar@{=>}"1";"2"^{\alpha} }$$
$$\xymatrix{
\cb(X,Y) \ar[r]^-{j\x 1}_{~}="1"  \ar@/_1pc/[dr]_{1}^(0.53){~}="2" & \cb(Y,Y)\x\cb(X,Y) \ar[d]^{M} \\
 & \cb(X,Y) 
\ar@{=>}"1";"2"^{\lambda} }$$
$$\xymatrix{
\cb(X,Y) \ar[r]^-{1\x j}_{~}="1" \ar@/_1pc/[dr]_{1}^(0.52){~}="2" & 
\cb(X,Y)\x\cb(X,X) \ar[d]^{M} \\
& \cb(X,Y) 
\ar@{=>}"2";"1"_{\rho} }$$
\end{itemize}
whose components take the form  $\alpha_{f,g,h}\colon(hg)f\to h(gf)$, $\lambda_f\colon 1f\to f$, and $\rho_f\colon f\to f1$, except that usually we omit the subscripts and simply write $\alpha$, $\lambda$, and $\rho$. These are required to satisfy five conditions, asserting the commutativity of all diagrams of the form
$$\xymatrix @C1pc @R1pc {
& (k(hg))f \ar[rr]^{\alpha_{f,hg,k}} & {}\ar@{}[dd]|{\fbox{1}} & k((hg)f) \ar[dr]^{k\alpha_{f,g,h}} & \\
((kh)g)f \ar[ur]^{\alpha_{g,h,k}f} \ar[drr]_{\alpha_{f,g,kh}} &&&& k(h(gf)) \\
&& (kh)(gf) \ar[urr]_{\alpha_{gf,h,k}} && }$$
$$\xymatrix{
(g1)f \ar[rr]^{\alpha_{f,1,g}} & \ar@{}[dd]|{\fbox{2}} & g(1f) \ar[dd]^{g\lambda_f} \\ \\
gf \ar[uu]^{\rho_g f} \ar@{=}[rr] && gf }
\xymatrix{
(1g)f \ar[rr]^{\alpha_{f,g,1}} \ar[dr]_{\lambda_g f} & \ar@{}[d]|{\fbox{3}} & 1(gf) \ar[dl]^{\lambda_{gf}} \\
& gf & }$$
$$\xymatrix{
&& gf \ar[dl]_{\rho_{gf}} \ar[dr]^{g\rho_f} \ar@{}[d]|{\fbox{4}} && 1 \ar@{=}[rr] \ar[dr]_{\rho_1} & \ar@{}[d]|{\fbox{5}} & 1 \\
&  (gf)1 \ar[rr]_{\alpha_{1,f,g}} && g(f1) && 11 \ar[ur]_{\lambda_1}  }$$
\begin{example}
In the usual way, we identify one-object skew bicategories with skew monoidal categories. 
\end{example}

\begin{example}\label{ex:bicategory}
  If the natural transformations $\alpha$, $\lambda$, and $\rho$ are all invertible, we recover the usual notion of bicategory, except that the usual definition includes only the first two axioms; but by adapting the argument of \cite{Kelly-unitaxioms} for monoidal categories, or applying the coherence theorem of \cite{MacLane-Pare}, one easily deduces that the other three axioms are a consequence of the first two.
\end{example}

\section{Skew warpings on skew bicategories} \label{sect:skew-warping}

In this section we make the basic definition which is a common generalization of skew warpings on skew monoidal categories, and pseudo mw-monads on bicategories.

\begin{definition}
  A {\em skew warping} on the skew bicategory \cb consists of:
\begin{itemize}
\item a function $D\colon\ob\cb\to\ob\cb$
\item functors $T\colon\cb(X,DY)\to\cb(DX,DY)$
\item 1-cells $K\colon X\to DX$ for each $X$
\item natural transformations
\end{itemize}
$$\xymatrix @C1.3pc {
\cb(Y,DZ)\x\cb(X,DY) \ar[rr]^{T\x T}_(0.49){~}="1" \ar[d]_{T\x 1} && 
\cb(DY,DZ)\x\cb(DX,DY) \ar[d]^{M} \\
\cb(DY,DZ)\x\cb(X,DY) \ar[r]_-{M} & 
\cb(X,DZ) \ar[r]_{T} \ar@{=>}"1"^{v} & \cb(DX,DZ) }$$
$$\xymatrix{
\cb(X,DY) \ar[r]^-{T\x K}_{~}="2" \ar@/_1pc/[dr]_1^(0.515){~}="1" & \cb(DX,DY)\x\cb(X,DX) \ar[d]^{M} \\
& \cb(X,DY) 
\ar@{=>}"1";"2"^{k} }\quad\quad 
\xymatrix{
1 \ar[r]^-{K}_(0.45){~}="1" \ar@/_1pc/[dr]_j^{~}="2"  & \cb(Y,DY) \ar[d]^{T} \\
& \cb(DY,DY) 
\ar@{=>}"1";"2"^{v_0} }$$
or, in terms of components
$$\xymatrix @R1pc {
T(Tg.f) \ar[r]^{v} & Tg.Tf \\
f \ar[r]^{k} & Tf.K \\
TK \ar[r]^{v_0} & 1_{DY} }$$
for $f\colon X\to DY$ and $g\colon Y\to DZ$.
\newline\noindent
These are required to satisfy the following five equations
$$\xymatrix{
T(T(Th.g).f) \ar[r]^{v} \ar[d]_{T(v.1)} & 
T(Th.g).Tf \ar[r]^{v.1}  \ar@{}[d]|{\fbox{1}} & 
(Th.Tg).Tf \ar[dr]^{\alpha} \\
T((Th.Tg).f) \ar[r]_{T\alpha} & 
T(Th.(Tg.f)) \ar[r]_{v} & Th.T(Tg.f) \ar[r]_{1.v} &
Th.(Tg.Tf) }$$
$$
\xymatrix{
T(Tf.K) \ar[r]^{v} \ar@{}[dr]|{\fbox{2}} & Tf.TK \ar[d]^{1.v_0} \\
Tf \ar[u]^{Tk} \ar[r]_{\rho} & Tf.1 }\quad\quad
\xymatrix{
T(TK.f) \ar[r]^{v} \ar[d]_{T(v_0.1)} & TK.Tf \ar[r]^{v_0.1} \ar@{}[d]|{\fbox{3}} & 
1.Tf \ar[d]^{\lambda} \\
T(1.f) \ar[rr]_{T\lambda} & & Tf }
$$
$$\xymatrix @R2pc {
T(Tg.f).K \ar[r]^{v.1} \ar@{}[dr]|{\fbox{4}} & (Tg.Tf).K \ar[d]^{\alpha} \\
Tg.f \ar[u]^{k} \ar[r]_{1.k} & Tg.(Tf.K) } \quad\quad\quad\quad
\xymatrix @R2pc {
TK.K \ar[r]^{v_0.1} \ar@{}[dr]|{\fbox{5}} & 1.K \ar[d]^{\lambda} \\
K \ar[u]^{k} \ar@{=}[r] & K }$$
for all $f\colon X\to DY$, $g\colon Y\to DZ$, and $h\colon Z\to DW$.
\end{definition}

\begin{example}
A skew warping on a skew monoidal category, in the sense of \cite{skew}, is literally the same as a skew warping on the corresponding one-object skew bicategory.
\end{example}

\begin{example}
  Any category can be seen as a skew bicategory with no non-identity 2-cells. A skew warping on a category is the same thing as an mw-monad on the category, and so amounts to an ordinary monad on the category. 
\end{example}

\begin{definition}
  A {\em warping} on a bicategory is a skew warping for which $v$, $k$, and $v_0$ are invertible. 
\end{definition}

\begin{example}
  A warping on a 2-category \cb is the same as a pseudo mw-monad (a no iteration pseudomonad in the language of \cite{MarmolejoWood-NoIterationPsm}). In more detail, $T$ is the functor $(~)^\bd$ of \cite{MarmolejoWood-NoIterationPsm}, while $K_X$ is the 1-cell $dX$. The 2-cells $\bd_A$, $\bd_f$, and $\bd_{f,h}$ of \cite{MarmolejoWood-NoIterationPsm} are the inverses of suitable components of our $v_0$, $k$, and $v$. Our five axioms are then conditions 8, 2, 3, 5, and 1 respectively of \cite{MarmolejoWood-NoIterationPsm}, while the remaining axioms 4, 6, and 7 of \cite{MarmolejoWood-NoIterationPsm} amount to naturality of $v$ and $k$. 
\end{example}

\section{The Kleisli construction for skew warpings}

We saw in Section~\ref{sect:monads} that the Kleisli category of a monad is easily constructed in terms of the corresponding mw-monad. We now describe an analogous construction for skew warpings; this is a straightforward generalization of \cite[Proposition~3.6]{skew}.

Given a skew warping, as in the previous section, there is a new skew bicategory $\cb_T$ with the same objects as \cb, and with hom-categories given by $\cb_T(X,Y)=\cb(X,DY)$. The composition functors are given by 
$$\xymatrix{
\cb(Y,DZ)\x \cb(X,DY) \ar[r]^-{T\x1} & \cb(DY,DZ)\x \cb(X,DY) \ar[r]^-{M} & \cb(X,DZ) }$$
so that the composite of $f\colon X\to DY$ and $g\colon Y\to DZ$ is $Tg\circ f\colon X\to DZ$. The identities are given by the $K\colon X\to DX$. The associativity maps have the form 
$$\xymatrix{
T(Th.g). f \ar[r]^{v. 1} & (Th. Tg). f \ar[r]^{\alpha} & Th.(Tg. f) }$$
and the identity maps have the form 
$$\xymatrix{
TK. f \ar[r]^-{v_0.1} & 1. f \ar[r]^{\lambda} & f & f \ar[r]^-{k} & Tf. K.}$$

\begin{remark}
We have numbered the axioms for skew bicategories and for skew warpings in such a way that to prove axiom $n$ for $\cb_T$ one needs only axiom $n$ for \cb and axiom $n$ for the skew warping.   
\end{remark}

\begin{proposition}
In the definition of a (skew) warping, if \cb is a bicategory and if $v$, $v_0$, and $k$ are invertible, then axioms 3, 4, and 5 follow from the first two axioms. 
\end{proposition}

\proof
Suppose that the first two axioms hold. Then we can still form the Kleisli construction $\cb_T$ as above, and the associativity and identity 2-cells will be invertible and satisfy axioms~1 and~2. Thus as explained in Example~\ref{ex:bicategory} this defines a bicategory, and the remaining (skew) bicategory axioms 3, 4, and 5 hold. Now axioms 4 and 5 for a skew warping are literally the same as axioms 4 and 5 for the skew bicategory $\cb_T$, while axiom 3 for a skew warping is a straightforward consequence of axiom 3 for the skew bicategory $\cb_T$.
\endproof

\begin{corollary}\label{cor:main}
A warping on a monoidal category, in the sense of \cite{BookerStreet-Tannaka}, is the same as a warping on the corresponding one-object bicategory, and so as a pseudomonad on the one-object bicategory.
\end{corollary}

\begin{corollary}
  Conditions 1, 3, and 5 in \cite[Definition~2.1]{MarmolejoWood-NoIterationPsm} follow from the other conditions. 
\end{corollary}

\section{Algebras}

We now generalize the definition of algebra given in \cite[Section~4]{MarmolejoWood-NoIterationPsm} to our setting. 

Let \cb be a skew bicategory, and consider a skew warping on \cb, as in Section~\ref{sect:skew-warping}.

\begin{definition}
  An {\em algebra} for the skew warping consists of an object $A\in\cb$ equipped with 
\begin{itemize}
\item a functor $E\colon\cb(X,A)\to\cb(DX,A)$ for each $X$
\item natural transformations
$$\xymatrix @C1.35pc {
\cb(Y,A)\x\cb(X,DY) \ar[r]^-{E\x1} \ar[dr]_{E\x T} & 
\cb(DY,A)\x \cb(X,DY) \ar[r]^-{M} \ar@{=>}[d]^{e} & 
\cb(X,A) \ar[d]^{E} \\
& \cb(DY,A)\x\cb(DX,DY) \ar[r]_-{M} & \cb(DX,A) }$$
$$\xymatrix{
\cb(Y,A) \ar[r]^-{E\x K}_(0.478){~}="1"  \ar@/_1pc/[dr]_{1}^{~}="2" & \cb(DY,A)\x \cb(Y,DY) \ar[d]^{M} 
\\
& \cb(Y,A) 
\ar@{=>}"2";"1"_{e_0} }$$
or in terms of components
$$\xymatrix @R1pc { 
E(Ea.x) \ar[r]^{e} & Ea.Tx \\
a \ar[r]^{e_0} & Ea.K }$$
where $a\colon Y\to A$ and $x\colon X\to DY$
\end{itemize}
subject to axioms asserting the commutativity of the following diagrams. 
$$\xymatrix{ 
E(E(Ea.x).y) \ar[r]^{e} \ar[d]_{E(e.1)} & 
E(Ea.x).Ty  \ar[r]^{e.1} & (Ea.Tx).Ty \ar[dr]^{\alpha} \\
E((Ea.Tx).y) \ar[r]_{E\alpha} & E(Ea.(Tx.y)) \ar[r]_{e} & 
Ea.T(Tx.y) \ar[r]_{1.v} & Ea.(Tx.Ty) }$$
$$\xymatrix{
E(Ea.K) \ar[r]^{e} & Ea.TK \ar[d]^{1.v_0} \\
Ea \ar[u]^{Ee_0} \ar[r]_{\rho} & Ea.1 }\quad\quad
\xymatrix{
E(Ea.x).K \ar[r]^{e.1} & (Ea.Tx).K \ar[d]^{\alpha} \\
Ea.x \ar[u]^{e_0} \ar[r]_{1.k} & Ea.(Tx.K) }$$
\end{definition}

\begin{example}
In the case of a warping on a 2-category, an algebra is the same as an algebra, in the sense of \cite[Section~4]{MarmolejoWood-NoIterationPsm}, for the corresponding pseudo mw-monad. Explicitly, in the definition of \cite{MarmolejoWood-NoIterationPsm} the functor $(~)^\ba$ is our $E$, while the 2-cells $\ba_h$ and $\ba_{g,h}$ are inverses of the components of our $e_0$ and $e$. Our three axioms are the axioms 6, 2, and 3 of \cite{MarmolejoWood-NoIterationPsm}; while the remaining three axioms of \cite{MarmolejoWood-NoIterationPsm} amount to naturality of $e$ and $e_0$.
\end{example}

\begin{proposition}
In the definition of algebra for a warping on a bicategory, the third axiom is a consequence of the other two.
\end{proposition}

\proof
We write as if the bicategory were strict. Consider the following diagram 
$$\xymatrix{
Ea.x \ar[r]^-{e_0} \ar[d]_{e_0} & E(Ea.x).K \ar[r]^-{e.1} \ar[d]_{e_0} & Ea.Tx.K \ar[d]^{e_0} \\
E(Ea.x).K \ar@{=}[d] \ar[r]^-{Ee_0.1} & E(E(Ea.x).K).K \ar[r]^-{E(e.1).1} \ar[d]_{e.1} & E(Ea.Tx.K).K \ar[dd]^{e.1} \\
E(Ea.x).K \ar[d]_{e.1} & E(Ea.x).TK.K \ar[d]_{e.1.1} \ar[l]_{1.v_0.1} \\
Ea.Tx.K & Ea.Tx.TK.K \ar[l]^-{1.1.v_0.1} & Ea.T(Tx.K).K \ar[l]^-{1.v.1} }$$
%
in which the large region in the bottom right corner commutes by the first equation (``the pentagon'') and the left central region commutes by the second equation (``the unit condition''), while all other regions commute by naturality. 

Since $e_0$ and $e.1$ are invertible we may cancel them, and conclude that the upper path in the diagram 
$$\xymatrix{
Ea.Tx.K \ar[r]^-{e_0} \ar@/_2pc/[rr]_{1.Tk.1} & 
E(Ea.Tx.K).K \ar[r]^{e.1} & 
Ea.T(Tx.K).K \ar[d]^{1.v.1} \\ 
&& Ea.Tx.TK.K \ar[r]^-{1.1.v_0.1} & Ea.Tx.K }$$
is the identity. But the lower path is also the identity, by the unit condition for the warping, so the two paths agree. Using invertibility of $v$ and $v_0$ we can cancel to obtain commutativity of the triangular region on the left. Thus the central triangular region in the diagram
$$\xymatrix{
&& E(Ea.x).K \ar[d]^(0.6){T(1.k).1} \ar[dr]^{e.1} \\
&&  E(Ea.Tx.K).K \ar[dr]^{e.1} & Ea.Tx.K \ar[d]^{1.Tk.1} \\
Ea.x \ar[r]_-{1.k} \ar[uurr]^{e_0}  & Ea.Tx.K \ar[ur]^{e_0} \ar[rr]_{1.Tk.1} && Ea.T(Tx.K).K 
}$$
also commutes, while the other regions commute by naturality. Cancelling $1.Tk.1$ gives the last equation. \endproof

\begin{corollary}
  The third axiom in the definition of \cite[Section~4]{MarmolejoWood-NoIterationPsm} is redundant.
\end{corollary}

\section{Formal mw-monads}

Monads can be defined in any bicategory \cite{bicategories} or indeed any skew bicategory, and the formal theory of monads in bicategories is well-understood \cite{ftm}. If $B$ is an object of a bicategory \ck, there is a monoidal structure on $\ck(B,B)$ with tensor product given by composition, and a monad in \ck on the object $B$ is a monoid in $\ck(B,B)$.

Here we sketch a setting for the formal theory of mw-monads. This has similarities with 
 \cite{AltenkirchChapmanUustalu}, although it differs both in the motivation and in the detail.

We write as if the bicategory \ck were strict. Let $i\dashv i^*$ be an adjunction in \ck, with $i\colon A\to B$. Then there is a skew monoidal structure on the hom-category $\ck(A,B)$, with tensor product $g\ox f$ given by $gi^*f$, and unit $i$. By associativity of \ck we have $(h\ox g)\ox f=hi^*gi^*f=h\ox(g\ox f)$, while $\lambda$ and $\rho$ are defined by 
$$\xymatrix{
i\ox f=ii^*f \ar[r]^-{\epsilon f} & f & f\ar[r]^-{f\eta} & fi^*i = f\ox i}$$
using the unit and counit of the adjunction $i\dashv  i^*$. 

A monoid in $\ck(A,B)$ consists of an arrow $d\colon A\to B$ equipped with maps $K\colon i\to d$ and $T\colon di^*d\to d$, satisfying the following three equations.
$$\xymatrix{
di^*di^*d \ar[r]^{T11} \ar[d]_{11T} & di^*d \ar[d]^{T} \\
di^*d \ar[r]_{T} & d }
\quad
\xymatrix{
ii^*d \ar[r]^{K11} \ar[dr]_{\epsilon 1} & di^*d \ar[d]^{T} \\
& d }
\quad
\xymatrix{
di^*i \ar[r]^{11K} & di^*d \ar[d]^{T} \\
d\ar[u]^{1\eta} \ar[r]_1 & d }$$


Composition with $i$ defines a functor $u=\ck(i,1)\colon\ck(B,B)\to\ck(A,B)$. For $f,g\colon B\to B$ we have 
$$\xymatrix{
u(g)\ox u(f) \ar@{=}[r] & gii^*fi \ar[r]^{1\epsilon11} & gfi \ar@{=}[r] & u(gf) }$$
while $u(1)=i$; this makes $u$ into a (normal) monoidal functor.
In particular, it sends monoids to monoids; that is, monads on $B$ to monoids in $\ck(A,B)$.

\begin{example}\label{ex:Prof}
Let \ck be the bicategory of profunctors. Recall that any functor $f\colon A\to B$ defines a profunctor $f_*\colon A\to B$ defined by $f_*(b,a)=B(b,fa)$, and that $f_*$ has a right adjoint $f^*$ defined by $f^*(a,b)=B(fa,b)$; we often write $f$ for $f_*$. Let $A$ be the discrete category on the same set of objects as $B$, and let $i$ be the inclusion. Then to give a functor  $d\colon A\to B$ and a 2-cell $K\colon i\to d$ in \ck is to give, for each object $x$ of $B$, an object $dx$ and a morphism $K\colon x\to dx$. To give $T\colon di^*d\to d$ is to give morphisms
$$ \int^{y\in A, a\in B} B(b,dy)\x B(iy,a)\x B(a,dx) \to B(b,dx)$$
natural in $b\in B$ and $x\in A$. Now naturality in $x$ and $y$ say nothing, since $A$ is discrete; while naturality in $a$ and $b$ reduce this, by Yoneda, to giving maps
$$T\colon B(iy,dx)\to B(dy,dx).$$
The three axioms for a monoid in $\ck(A,B)$ are exactly the three axioms for an mw-monad. Thus a functor $A\to B$ is a monoid in $\ck(A,B)$ precisely when it is an mw-monad. Moreover, given this identification, the monoidal functor $u=\ck(i,B)$ sends a monad on $B$ to the corresponding mw-monad.
\end{example}

Motivated by this example, we consider monoids in $\ck(A,B)$ as our formal notion of mw-monad; of course monoids in $\ck(B,B)$ are our formal notion of monad. (This notion of mw-monad depends on $A$ and $i$, somewhat as in the treatment of  \cite{AltenkirchChapmanUustalu}.) 

In order to compare monads with mw-monads in this formal context, we should therefore compare monoids in $\ck(B,B)$ with monoids in $\ck(A,B)$. In the following section we propose a more general setting in which to perform this comparison.

\section{Normalization}

In this section we show that, under mild conditions, a skew monoidal category \cc can be replaced by a {\em right normal} skew monoidal category, meaning one for which the right unit map $\rho$ is invertible. Furthermore, the two skew monoidal categories have equivalent categories of monoids. We use this to complete the comparison between monads and mw-monads begun in the previous section.

Let \cc be a skew monoidal category with tensor $\ox$ and unit $I$; we shall often write $XY$ for $X\ox Y$. Suppose that \cc has reflexive coequalizers, and that these are preserved by tensoring on the right. The functor $-\ox I\colon\cc\to\cc$ given by tensoring on the right with the unit $I$ underlies a monad (see \cite{Szlachanyi-skew}) with the maps 
$$\xymatrix{
(X\ox I)\ox I \ar[r]^{\alpha} & X\ox (I\ox I) \ar[r]^-{1\ox\lambda} & X\ox I &
X \ar[r]^{\rho} & X\ox I }$$
defining the components of the multiplication and unit.  Write $\cc^I$ for the category of algebras for the monad; we call its objects $I$-modules. This has reflexive coequalizers, formed as in \cc.

If $(Y,y)$ is an $I$-module, and $X$ an arbitrary object of \cc, then $X\ox Y$ becomes an $I$-module via the action
$$\xymatrix{
(XY)I \ar[r]^{\alpha} & X(YI) \ar[r]^{1y} & XY, }$$
with associativity and unit axioms proved using the following diagrams.
$$\xymatrix{
((XY)I)I \ar[r]^{\alpha1} \ar[dd]_{\alpha} & (X(YI))I \ar[r]^{(1y)1} \ar[d]^{\alpha} & (XY)I \ar[d]^{\alpha}  \\
& X((YI)I) \ar[r]^{1(y1)} \ar[d]^{1\alpha} & X(YI) \ar[d]^{1y} \\
(XY)(II) \ar[r]^{\alpha} \ar[d]_{(11)\lambda} & X(Y(II)) \ar[d]_{1(1\lambda)} & XY \ar@{=}[d] \\
(XY)I \ar[r]_{\alpha} & X(YI) \ar[r]_{1y} & XY 
}
\quad
\xymatrix{
XY \ar[r]^{\rho} \ar[dr]_{1\rho} \ar@/_1pc/[ddr]_1 & (XY)I \ar[d]^{\alpha} \\
& X(YI) \ar[d]^{1y} \\
& XY }
$$

Given $I$-modules $(X,x)$ and $(Y,y)$, we may form the reflexive coequalizer
\begin{equation}
  \label{eq:rc}
\xymatrix @R1pc {
(XI)Y \ar[rr]^{x1} \ar[dr]_{\alpha} && XY \ar@/_2pc/[ll]_{\rho1} \ar[r]^{q} & X\wedge Y \\
& X(IY) \ar[ur]_{1\lambda} }
\end{equation}
in \cc, and this lifts to a coequalizer in the category of $I$-modules, whose object-part involves an action  $c\colon (X\wedge Y)I\to X\wedge Y$. This defines a functor $\wedge\colon\cc^I\x \cc^I\to\cc^I$. By commutativity of the diagram
$$\xymatrix @R1.5pc {
& ((XI)Y)Z \ar[r]^{\alpha1} \ar[dl]_{(x1)1} \ar[dd]^{\alpha} & (X(IY))Z \ar[dr]^{(1\lambda)1} \ar[d]^{\alpha} \\
(XY)Z \ar[dd]_{\alpha} && X((IY)Z) \ar[d]_{1\alpha} \ar[dr]^{1(\lambda1)} & (XY)Z \ar[d]^{\alpha} \\
& (XI)(YZ) \ar[dl]_{x(11)} \ar[d]^{(11)q} \ar[r]^{\alpha} & X(I(YZ)) \ar[d]^{1(1q)} \ar[r]^{1\lambda} & X(YZ) \ar[d]^{1q} \\
X(YZ) \ar[d]_{1q} & (XI)(Y\wedge Z) \ar[dl]_{x1} \ar[r]^{\alpha} & X(I(Y\wedge Z)) \ar[r]^{1\lambda} & X(Y\wedge Z)  \ar[dl]^q \\
X(Y\wedge Z) \ar[rr]^{q} && X\wedge(Y\wedge Z) }$$
%
%
there is a unique induced $\alpha_1\colon (X\wedge Y)Z\to X\wedge(Y\wedge Z)$ whose composite with $q1\colon(XY)Z\to (X\wedge Y)Z$ is $q.1q.\alpha$.  The various regions of the diagram
$$\xymatrix @C1.5pc {
&& ((X\wedge Y)I)Z \ar[r]^{c1} & (X\wedge Y)Z \ar@/^4pc/[ddddd]^{\alpha_1} \\
& ((XY)I)Z \ar[ur]^{(q1)1} \ar[r]_{\alpha1} \ar[dd]^{\alpha} \ar[dl]_{(q1)1} & 
(X(YI))Z \ar[r]_{(1y)1} \ar[d]^{\alpha} & (XY)Z \ar[u]^{q1} \ar[d]_{\alpha} &  \\
((X\wedge Y)I)Z \ar[dd]_{\alpha} && X((YI)Z) \ar[r]_{1(y1)} \ar[d]^{1\alpha} & 
X(YZ) \ar[dd]_{1q} \\
& (XY)(IZ) \ar[dl]_{q(11)} \ar[r]^{\alpha} \ar[d]^{(11)\lambda} & X(Y(IZ)) \ar[d]^{1(1\lambda)} \\
(X\wedge Y)(IZ) \ar[d]_{1\lambda} & (XY)Z \ar[r]^{\alpha} \ar[dl]^{q1} & X(YZ) \ar[r]^{1q} & X(Y\wedge Z) \ar[d]_{q}    \\
(X\wedge Y)Z \ar@/_1pc/[rrr]_{\alpha_1} &&&  X\wedge(Y\wedge Z)
}$$
are easily seen to commute, thus the exterior does so. Cancel the epimorphism $(q1)1$, and deduce the commutativity of the diagram which guarantees that $\alpha_1$ factorizes uniquely through $q\colon(X\wedge Y)Z\to (X\wedge Y)\wedge Z$ to give a morphism $\alpha'\colon (X\wedge Y)\wedge Z\to X\wedge (Y\wedge Z)$ making the triangle in the diagram
$$\xymatrix{
(XY)Z \ar[d]_{\alpha} \ar[r]^{q1} & (X\wedge Y)Z \ar[r]^{q} \ar[dr]_{\alpha_1} & (X\wedge Y)\wedge Z \ar[d]^{\alpha'} \\
X(YZ) \ar[r]_{1q} & X(Y\wedge Z) \ar[r]_q & X\wedge(Y\wedge Z) }$$
commute. The larger region on the left commutes by definition of $\alpha_1$, and so the exterior commutes.

The resulting $\alpha'$ is clearly natural, and commutativity of the pentagon for $\alpha$ implies commutativity of the pentagon for $\alpha'$.

Commutativity of the diagrams
$$\xymatrix{
(II)I \ar[r]^{\alpha} \ar[d]_{\lambda 1} & I(II) \ar[r]^{1\lambda} \ar[dl]^{\lambda} & II \ar[d]^\lambda \\
II \ar[rr]_\lambda  && I 
}
\quad
\xymatrix{
I \ar[r]^{\rho} \ar[dr]_1 & II \ar[d]^\lambda \\ & I }$$
shows that $\lambda\colon II\to I$ makes $I$ into an $I$-module. 

Commutativity of 
$$\xymatrix{
(II)X \ar[d]_{\alpha} \ar[r]^{\lambda 1} & IX \ar[dd]^\lambda \\
I(IX) \ar[ur]_\lambda \ar[d]_{1\lambda} \\
IX \ar[r]_{\lambda} & X }$$
shows that $\lambda\colon IX\to X$ factorizes uniquely through $q\colon IX\to I\wedge X$ to give a map $\lambda'\colon I\wedge X\to X$.

On the other hand, the diagram
$$\xymatrix{
(XI)I \ar[rr]_{x1} \ar[dr]_{\alpha} && XI \ar[r]_{x} \ar@/_1pc/@{.>}[ll]_{\rho} & X \ar@/_1pc/@{.>}[l]_{\rho} \\
& X(II) \ar[ur]_{1\lambda} 
}$$
is a split coequalizer in \cc, and the solid part is a fork in $\cc^I$, thus is a coequalizer in $\cc^I$, and so  exhibits $X$ itself as $X\wedge I$. Rather than  identify $X\wedge I$ with $X$, though, we  let $\rho'$ be the composite 
$$\xymatrix{
X \ar[r]^{\rho} & XI \ar[r]^{q} & X\wedge I}$$
and note that this is invertible. 

We now show that $\alpha'$, $\rho'$, and $\lambda'$ make $\cc^I$ into a skew monoidal category. We have already observed that  the pentagon commutes, so we turn to the four remaining axioms.

Compatibility of $\alpha'$ and $\rho'$ follows from the corresponding condition for $\alpha$ and $\rho$, and commutativity of the diagrams
$$\xymatrix{
XY \ar[d]_q \ar[r]^\rho & (XY)I \ar[d]_{q1} \ar[r]^{\alpha} & X(YI) \ar[d]^{1q} \\
X\wedge Y \ar[r]^\rho \ar[dr]_{\rho'} & (X\wedge Y)I \ar[d]_{q} & X(Y\wedge I) \ar[d]^q \\
& (X\wedge Y)\wedge I \ar[r]^{\alpha'}  & X\wedge(Y\wedge I) }
~~~~
\xymatrix{
XY \ar[d]_q \ar[r]^{1\rho} \ar[dr]_{1\rho'} & X(YI) \ar[d]^{1q} \\
X\wedge Y \ar[dr]_{1\wedge\rho'} & X(Y\wedge I) \ar[d]^q \\
& X\wedge(Y\wedge I) }$$

Compatibility of $\alpha'$ and $\lambda'$ follows from the corresponding condition for $\alpha$ and $\lambda$, and commutativity of the diagrams
$$\xymatrix{
(IX)Y \ar[d]_{q1} \ar[r]^\alpha & I(XY) \ar[d]^{1q} \ar[dr]^{\lambda} \\
(I\wedge X)Y \ar[d]_q & I(X\wedge Y) \ar[d]^{q} \ar[dr]^{\lambda} & XY \ar[d]^q \\
(I\wedge X)\wedge Y \ar[r]_{\alpha'} & I\wedge (X\wedge Y) \ar[r]_{\lambda'} & X\wedge Y }~~~ 
\xymatrix{
(IX)Y \ar[d]_{q1} \ar[dr]^{\lambda1} \\
(I\wedge X)Y \ar[d]_q \ar[r]_-{\lambda'1} & XY \ar[d]^q \\
(I\wedge X)\wedge Y \ar[r]_-{\lambda'\wedge1} & X\wedge Y }$$

For the triple compatibility condition, observe that the diagram 
$$\xymatrix{
XY \ar[r]^q \ar[dr]^{\rho'1} \ar[d]_{\rho1} & X\wedge Y \ar[dr]^{\rho'\wedge1} \\
(XI)Y \ar[r]_{q1} \ar[d]_{\alpha} & (X\wedge I)Y \ar[r]_{q} & 
(X\wedge I)\wedge Y \ar[d]^{\alpha'} \\
X(IY) \ar[r]^{1q} \ar[dr]_{1\lambda} & X(I\wedge Y) \ar[d]^{1\lambda'} \ar[r]^{q} & X\wedge (I\wedge Y) \ar[d]^{1\wedge\lambda'} \\
& XY \ar[r]_{q} & X\wedge Y }$$
commutes and that $q$ is epi; then the axiom for $\cc^I$ follows from that for \cc. 

Finally compatibility of $\lambda'$ and $\rho'$ follows from commutativity of 
$$\xymatrix{
& I\wedge I \ar[dr]^{\lambda'} \\
I \ar[ur]^{\rho'} \ar[r]^{\rho} \ar@/_1pc/[rr]_{1} & 
II \ar[u]_q \ar[r]^{\lambda} & I. }$$

This now proves that $\cc^I$ is skew monoidal; indeed it is a right normal skew monoidal category, in the sense that $\rho$ is invertible.  The forgetful functor $U\colon\cc^I\to\cc$ is a monoidal functor, with $U_2\colon U(X,x)\ox U(Y,y) \to U((X,x)\wedge(Y,y))$ given by the quotient map $q\colon XY\to X\wedge Y$, and $U_0\colon I\to U(I,\lambda)$ the identity.

This process is universal, in the sense that if \cd is any right normal skew monoidal category and $M\colon\cd\to\cc$ a monoidal functor, then $M$ factorizes uniquely through $U$ as a skew monoidal functor $N\colon\cd\to\cc^I$. For each object $X\in\cd$, we have an $I$-module structure on $MX$, given by 
$$\xymatrix @C3pc {
MX\ox I \ar[r]^-{1\ox M_0} & MX\ox MI \ar[r]^{M_2} & 
M(X\ox I) \ar[r]^-{M(\rho^{-1})} & MX }$$
and this is natural in $X$, so that $M$ does lift to functor $N\colon\cd\to\cc^I$ with $UN=M$. 

Furthermore, by commutativity of 
$$\xymatrix{
(MX.I).MY \ar[r]^{(1M_0)1} \ar[d]_\alpha & (MX.MI).MY \ar[r]^{M_2.1} \ar[d]_\alpha & 
M(XI).MY \ar[r]^{M\rho^{-1}.1} \ar[d]^{M_2} & MX.MY \ar[d]^{M_2} \\
MX.(I.MY) \ar[r]^{1(M_01)} \ar[d]_{1\lambda} & MX.(MI.MY) \ar[d]^{1.M_2} & M((XI)Y) \ar[r]^{M(\rho^{-1}1)} \ar[d]^{M\alpha} & M(XY) \ar@{=}[d] \\
MX.MY \ar@/_2pc/[rrr]_{M_2} & MX.M(IY) \ar[r]^{M_2} \ar[l]_{1.M\lambda} & M(X(IY)) \ar[r]^{M(1\lambda)} & M(XY) 
}$$
we see that $M_2$ passes to the quotient to give a map $N_2\colon NX\wedge NY\to N(XY)$; while $M_0$ underlies a map $N_0\colon I\to NI$.

Since monoids in a skew monoidal category are just monoidal functors out of the terminal skew monoidal category, and this terminal skew monoidal category is right normal (in fact monoidal), it follows that the monoids in $\cc^I$ are the same as the monoids in \cc.

We summarize this as follows:

\begin{theorem}
  Let \cc be a skew monoidal category, and suppose that \cc has coequalizers of reflexive pairs of the form \eqref{eq:rc}, and that these are preserved by tensoring on the right. Then the category $\cc^I$ of right $I$-modules is a right normal skew monoidal category, and the forgetful functor $U\colon\cc^I\to\cc$ is a normal monoidal functor. Furthermore, it is universal, in the sense that for any right normal skew monoidal category \cd, composition with $U$ induces an equivalence between the category of monoidal functors from \cd to $\cc^I$ and the category of monoidal functors from \cd to \cc.
\end{theorem}

We call $\cc^I$ the {\em right normalization} of the skew monoidal category \cc.

The next result is our promised formal approach to the comparison of monads and mw-monads.

\begin{theorem}
  Let  $i\colon A\to B$ be a morphism in a bicategory \ck, and suppose that $i$ has a right adjoint $i\dashv i^*$ and is opmonadic (of Kleisli type). Suppose further that for any $h\colon B\to B$, the functor $\ck(h,B)\colon\ck(B,B)\to\ck(B,B)$ preserves any existing coequalizers of reflexive pairs. Then the skew monoidal category $\ck(A,B)$ satisfies the conditions of the previous theorem, and the right normalization $\ck(A,B)^I$ is given by $\ck(B,B)$. Thus monoids in $\ck(A,B)$ are equivalent to monoids in $\ck(B,B)$. 
\end{theorem}

\proof
The adjunction $i\dashv i^*$ induces an adjunction $\ck(i^*,B)\dashv\ck(i,B)$, which in turn induces a monad on $\ck(A,B)$, and this monad is precisely that given by tensoring on the right with the unit $i$ of $\ck(A,B)$. Since $i$ is opmonadic, $\ck(i,B)$ is monadic, and so $\ck(B,B)$ is equivalent to the category of $I$-modules. 

Using again the fact that $i$ is opmonadic, the diagram 
$$\xymatrix{
gii^*ii^* \ar@<1ex>[r]^{g\epsilon ii^*} \ar@<-1ex>[r]_{gii^*\epsilon} & gii^* \ar[r]^{g\epsilon} & g }$$
is a coequalizer in $\ck(B,B)$, and now composing on the right with $fi$, we see that the required coequalizers \eqref{eq:rc} exist, with $gi\wedge fi=(gf)i$.

Thus the normalization does exist, and since $u\colon\ck(B,B)\to\ck(A,B)$ is a monoidal functor with right normal domain (in fact monoidal domain), we have the comparison $v\colon\ck(B,B)\to\ck(A,B)^I$. From the construction of $v$ it is clear that this is a monoidal equivalence. 
\endproof

\begin{example}
Consider the case of Example~\ref{ex:Prof}, where \ck is the bicategory of profunctors, and $i\colon A\to B$ is the inclusion of the discrete category $A$ on the same set of objects as $B$. Since $i$ is the identity on objects it is indeed opmonadic, while $\ck(A,B)$ is cocomplete with colimits preserved by tensoring on either side, thus the conditions of the theorem hold. We recover the correspondence between monads and mw-monads by observing that a profunctor $g\colon B\to B$ is a functor if and only if the composite $gi\colon A\to B$ is one. 
\end{example}


\providecommand{\bysame}{\leavevmode\hbox to3em{\hrulefill}\thinspace}
\providecommand{\MR}{\relax\ifhmode\unskip\space\fi MR }
\providecommand{\MRhref}[2]{%
  \href{http://www.ams.org/mathscinet-getitem?mr=#1}{#2}
}
\providecommand{\href}[2]{#2}

\end{document}